\documentclass[11pt,reqno]{amsart}
\usepackage[margin=1in]{geometry}
\usepackage{amsmath,amssymb,amsthm,mathtools}
\usepackage[T1]{fontenc}
\usepackage{lmodern}
\usepackage{enumitem}
\usepackage[colorlinks=true,linkcolor=blue,citecolor=blue,urlcolor=blue]{hyperref}
\hypersetup{pdftitle={When are noetherian Hopf
algebras affine?},
pdfauthor={Hongdi Huang, Xin Tang, Xingting Wang, James J. Zhang},
pdfkeywords={Hopf algebra, noetherian ring, affine algebra, free subalgebra, polynomial identity, Gelfand-Kirillov dimension}}
\usepackage{microtype}

\newtheorem{theorem}{Theorem}[section]
\newtheorem{proposition}[theorem]{Proposition}
\newtheorem{lemma}[theorem]{Lemma}

\theoremstyle{definition}
\newtheorem{definition}[theorem]{Definition}
\newtheorem{corollary}[theorem]{Corollary}
\newtheorem{question}[theorem]{Question}
\newtheorem{example}[theorem]{Example}

\numberwithin{equation}{section}

\newcommand{\kk}{\Bbbk}
\newcommand{\eps}{\epsilon}
\DeclareMathOperator{\GKdim}{GKdim}
\DeclareMathOperator{\injdim}{injdim}

\newcommand{\Ql}{Q_{\mathrm{l}}}

\title{When are noetherian Hopf
algebras affine?}

\author{Hongdi Huang}
\address{Department of Mathematics, Shanghai University, Shanghai 200444, China. Newtouch Center for Mathematics of Shanghai University,
Shanghai 200444, China}
\email{hdhuang@shu.edu.cn}
\author{Xin Tang}
\address{Department of Mathematics \& Computer Science,
Fayetteville State University, Fayetteville, NC 28301, USA}
\email{xtang@uncfsu.edu}
\author{Xingting Wang}
\address{Department of Mathematics, Louisiana State University,
Baton Rouge, Louisiana 70803, USA}
\email{xingtingwang@math.lsu.edu}
\author{James J. Zhang}
\address{Department of Mathematics, Box 354350,
University of Washington, Seattle, Washington 98195, USA}
\email{zhang@math.washington.edu}

\begin{document}
\begin{abstract}
We show that a noetherian Hopf algebra over 
an arbitrary base field is finitely 
generated as an algebra when it contains 
no free subalgebra of rank two. In 
particular, noetherian Hopf algebras 
that are PI or have finite
Gelfand--Kirillov dimension are 
affine, which answers a question of 
Wu and the fourth-named author. 
Consequently, noetherian PI Hopf 
algebras are AS--Gorenstein, 
Auslander--Gorenstein, GK--Cohen--Macaulay 
and residually finite-dimensional. 
\end{abstract}
\maketitle

\section*{Introduction}
\label{xxsec0}
\label{sec:intro}

Throughout, $\kk$ is a field of arbitrary 
characteristic. All algebras are assmed to be unital 
and associative over $\kk$ unless 
otherwise stated, and unadorned tensor 
products are over $\kk$. An algebra is 
called \emph{affine} if it is finitely
generated as an algebra over $\kk$. We 
refer to \cite{Mo93} for basic
definitions about Hopf algebras.

\subsection{Motivation and background}
\label{xxsec0.1}
There has been a wealth of examples of noetherian
Hopf algebras, including group algebras 
of polycyclic-by-finite groups
\cite{Ha54}, enveloping algebras of 
finite-dimensional Lie algebras and
coordinate rings of affine algebraic 
groups. The list is further expanded after the 
``birth'' of quantum groups. In particular, quantized
enveloping algebras of finite-dimensional semisimple Lie algebras and quantized 
coordinate rings for semisimple algebraic groups are noetherian Hopf 
algebras, and they satisfy a polynomial identity 
(PI) in the roots of unity case \cite[Part~III]{BG02}. Brown and 
Goodearl initiated the systematic study of 
noetherian PI Hopf algebras \cite{BG97}; and at the 1997 Seattle conference Brown 
proposed a series of problems and 
conjectures \cite{Br98}. The list of open questions has 
kept growing over the years; we refer to the
surveys \cite{Br07,Go13,BGi14,BZ21,And25,Br26} 
for the development and status of this subject.

In this note we study the following question, which was raised earlier 
by Wu and the fourth-named author \cite[Question 5.1]{WZ03}.

\begin{question}
\label{xxque0.1}
Is every noetherian \textup{(}resp. noetherian and PI\textup{)} 
Hopf algebra affine?
\end{question}

The same question appears in 
\cite[Question C(ii)]{Br07},
\cite[Question 2.4(a)]{Go13}, 
\cite[Question D]{BGi14},
\cite[Question 6]{And25} and 
\cite[Question 2.1]{Br26}. Much of the
structure theory for noetherian Hopf 
algebras has been developed under the
additional assumption that the Hopf 
algebra is affine. For example, Wu
and the fourth-named author proved 
that every affine noetherian PI Hopf
algebra is AS-Gorenstein, 
Auslander--Gorenstein and 
GK--Cohen--Macaulay
\cite[Theorems 0.1 and 0.2]{WZ03}, 
answering a question of Brown
\cite[Question A]{Br98}. A positive 
answer to Question \ref{xxque0.1} would remove
the affineness hypothesis from such results.

The answer to Question \ref{xxque0.1} is known to be positive in 
several cases. Molnar proved that
commutative 
noetherian Hopf algebras are affine
\cite{Mo75}, and artinian Hopf algebras 
are proved to be finite dimensional
\cite{LZ07}. Left noetherian pointed 
Hopf algebras are affine
\cite[Theorem 1.5]{GZ17}, \cite{JZ25}; 
and over an algebraically closed
field, a commutative-by-finite Hopf 
algebra is noetherian if and only if
it is affine \cite[Theorem 2.3]{BC21}. 
Jia and Y. Zhang proved
affineness for certain noetherian 
graded Hopf algebras \cite{JZ26a}, and
showed that a one-sided noetherian 
Hopf algebra is affine if and only if
its Hopf coradical is affine \cite{JZ26b}. 
Recall that an algebra is
\emph{residually finite-dimensional} if the intersection of all its ideals
with finite codimension is zero. Skryabin proved that a residually
finite-dimensional noetherian PI Hopf algebra is finitely generated as a Hopf
algebra \cite[Proposition 4.8]{Sk21b}.

Beyond these cases, the situation remains muddy. Note that Molnar's proof employs two
essential ingredients: the antipode of a commutative Hopf algebra has finite order, and
Hopf subalgebras correspond to Hopf ideals. Neither is available for
general noetherian PI Hopf algebras \cite[\S5]{WZ03}. Nor does small growth
force affineness: Goodearl and the fourth-named author constructed Hopf
domains of Gelfand--Kirillov dimension one and two, which are neither
affine nor noetherian \cite{GZ17}. Even for PI Hopf algebras, the question has remained open
\cite[Remark 3.8(i)]{Br26}, \cite{JZ26b}.

\subsection{Main results}
\label{xxsec0.2}
Skryabin conjectured that every noetherian Hopf algebra has bijective
antipode \cite{Sk06}. This conjecture has recently been settled by Sciandra \cite[Theorem 6.1]{Sc26}. As a consequence, a Hopf algebra is left
noetherian if and only if it is right noetherian
\cite[Corollary 6.2]{Sc26}, so we simply say \emph{noetherian}. Two elements $x$ and $y$ of an algebra form a \emph{free
pair} if the subalgebra they generate is free on $x$ and $y$. Our main
result is the following.

\begin{theorem}
\label{xxthm0.2}
Let $H$ be a noetherian Hopf algebra. If $H$ does not contain a free pair, then $H$ is affine.
\end{theorem}

The free algebra $\kk\langle x,y\rangle$ does not satisfy any nontrivial polynomial identity
and has exponential growth. Recall that an algebra is \emph{locally PI}
if every affine subalgebra of it is PI. So Theorem~\ref{xxthm0.2} applies
in the following cases.

\begin{corollary}
\label{xxcor0.3}
Let $H$ be a noetherian Hopf algebra. Then $H$ is affine if either $H$ is
locally PI or every affine subalgebra of $H$ has subexponential
growth. In particular, $H$ is affine if it either satisfies a polynomial
identity or has finite Gelfand--Kirillov dimension.
\end{corollary}

Corollary~\ref{xxcor0.3} answers \cite[Question 5.1]{WZ03} and
\cite[Question C(ii)]{Br07} for PI Hopf algebras and Hopf algebras of
finite Gelfand--Kirillov dimension respectively. It also partially answers \cite[Question~1.4(4)]{GZ17},
which asks whether affineness and noetherianity are equivalent
for Hopf domains of finite Gelfand--Kirillov dimension,
by establishing the implication from noetherianity to affineness.

\subsection{Strategy of the proof}
\label{xxsec0.3}
To aid the proof, we introduce the notions of \emph{epimorphic} subalgebras
(Definition~\ref{xxdef1.1}) and \emph{\textup{(}regular\textup{)} pseudo-birational}
subalgebras
(Definition~\ref{xxdef1.3}). Let $N$ be the prime radical of $H$. The
proof of Theorem~\ref{xxthm0.2} goes in four steps:
\begin{equation}
\label{E0.3.1}\tag{E0.3.1}
\begin{array}{r@{\;\;}c@{\;\;}l}
\text{$H$ is noetherian}
&\xRightarrow{\ \text{Lemma~\ref{xxlem1.2}(2)}\ }
&\text{$\exists$ affine epimorphic subalgebra of $H$}\\[3pt]
&\xRightarrow{\ \text{Lemmas~\ref{xxlem1.5}(1), \ref{xxlem2.5}}\ }
&\text{$\exists$ affine regular pseudo-birational subalgebra of $H/N$}\\[3pt]
&\xRightarrow{\ \text{Lemma~\ref{xxlem1.5}(3)}\ }
&\text{$\exists$ affine pseudo-birational subalgebra of $H$}\\[3pt]
&\xRightarrow{\ \text{Lemma~\ref{xxlem3.1}(3)}\ }
&\text{$H$ is affine.}
\end{array}
\end{equation}

Step 1 only uses the fact that the augmentation ideal $H^+$ is a finitely generated left ideal of $H$. Step 2
is the only place where the absence of free pair is used. Through Jategaonkar's theorem \cite{Ja69}, it furnishes a suitable affine subalgebra
$B$ of $H/N$ with a semisimple artinian classical left quotient ring, and the
epimorphic property then forces $\Ql(B)=\Ql(H/N)$. Step 3 lifts common
left denominators through the nilpotent ideal $N$. Since these denominators are multiplied together, they must be regular, not merely nonzero. Step 4 does need the 
antipode to be bijective, which is Sciandra's theorem
\cite{Sc26}.

We carry out these four steps in Sections~\ref{xxsec1} through \ref{xxsec3}. Section~\ref{xxsec4} contains further comments, examples and questions.


\section{Epimorphic and (regular) pseudo-birational subalgebras}
\label{xxsec1}

\begin{definition}
\label{xxdef1.1}
A subalgebra $A$ of an algebra $\Lambda$ is called {\it epimorphic} if the
multiplication map $\Lambda \otimes_A \Lambda \to \Lambda$ is bijective. 
\end{definition}

The above definition of an epimorphic subalgebra is equivalent to the condition that the inclusion $A\rightarrow \Lambda$ is an epimorphism in the category of rings \cite[Proposition~XI.1.2]{Ste75}. For a Hopf algebra $H$, we write $\Delta$, $\eps$ and $S$ for the coproduct, counit and antipode of $H$ respectively. We adopt Sweedler notation and set $H^+=\ker\eps$ (even when $H$ is only a bialgebra). A \emph{subbialgebra} of a Hopf algebra $H$ is a subalgebra which is also a subcoalgebra; it need not be stable under $S$. 

For any subalgebra $A$ of $H$, we denote by $A^{+}\colon =H^{+}\cap A$. The first step in \eqref{E0.3.1}, namely,
$${\text{$H$ being  noetherian}}
\Rightarrow
{\text{$\exists$ affine epimorphic subalgebra}}$$
is part (2) of 
the following lemma.

\begin{lemma}
\label{xxlem1.2}
Let $A$ be a subbialgebra of a Hopf algebra $H$.
\begin{enumerate}
\item[{\rm{(1)}}]
The following are equivalent.
\begin{enumerate}
\item[(1a)]
$A$ is an epimorphic subbialgebra.
\item[(1b)]
$HA^{+}=H^{+}$.
\item[(1c)]
$A^{+}H=H^{+}$.
\end{enumerate}
\item[{\rm{(2)}}]
The following are equivalent.
\begin{enumerate}
\item[(2a)]
$H^{+}$ is finitely generated as a left ideal of $H$. 
\item[(2b)]
$H$ contains an affine epimorphic subbialgebra.
\end{enumerate}
As a consequence, if $H$ is left noetherian, then it contains an affine epimorphic subbialgebra.

\end{enumerate}
\end{lemma}


\begin{proof}
(1) Since $A$ is a subcoalgebra of $H$, we have that $h_1a_1-h_1\eps(a_1)\in HA^+$ for $a\in A$. Thus $\Phi(h\otimes_A x)=\sum\overline{h_1}\otimes h_2x$ defines a well-defined
map $\Phi\colon H\otimes_AH\to(H/HA^+)\otimes H$. 
Note that the map 
$\Phi$ is similar to the map 
$\beta$ defined 
in \cite[p.35, l.7]{Mo93}.
The map $\Psi\colon (H/HA^+)\otimes H\to H\otimes_AH$ with $\Psi(\overline h\otimes x)=\sum h_1\otimes_AS(h_2)x$, is also well defined,
because for any $a\in A^+$, one has
\[
\sum h_1a_1\otimes_AS(a_2)S(h_2)x
=\sum h_1\otimes_Aa_1S(a_2)S(h_2)x
=\sum h_1\otimes_A\epsilon(a)S(h_2)x=0.
\]
It is straightforward to check that $\Phi$ and $\Psi$ are inverses of each other. As a result, both $\Phi$ and $\Psi$ are bijective. 

(1b) $\Rightarrow$ (1a) Since
$HA^+=H^+$, one can use the counit to identify $H/HA^+$ with $\kk$. As a result, $\Phi$ can be identified with the multiplication: $H\otimes_{A}H\to H$, which is thus bijective. 

(1a) $\Rightarrow$ (1b) Let $m: H\otimes_A H\to H$ be the multiplication. Then 
$m\circ \Psi$ maps the element $\bar{h}\otimes x$ to $\epsilon(h)x$. Since $A$ is epimorphic,
$m\circ \Psi$ is an isomorphism, which implies
that $H/HA^{+}=\Bbbk$. Therefore, $H^{+}=HA^{+}$.

Similarly, (1a) $\Leftrightarrow$ (1c).

(2) Suppose $H^{+}$ is generated by a finite set
$X$ as a left ideal of $H$. Let $C$ be a  
finite dimensional coalgebra containing 
$X$, and denote by $A$ the subalgebra of $H$ generated by $C$. Then $A$ is affine and $HA^{+}=H^{+}$. By part (1), $A$ is an epimorphic subbialgebra of $H$. Conversely, assume that
$A$ is an affine 
epimorphic subbialgebra of $H$ 
(or equivalently $H^{+}=HA^{+}$). Let $X$ be a finite set of algebra generators for $A$. By replacing $x\in X$ by $x-\epsilon(x)$,
we may assume $X\subseteq H^{+}$. Then 
$H^{+}=HA^{+}=\sum_n HX^{n}=HX$, which is finitely
generated as a left ideal of $H$.
\end{proof}

An element of a ring is called \emph{regular} if it is neither a left nor a
right zero divisor. 

\begin{definition}
\label{xxdef1.3}
Let $\Lambda$ be an algebra and $A$ be a subalgebra of $\Lambda$.
\begin{enumerate}
\item[(1)]
We say $A$ is {\it pseudo-birational}
in $\Lambda$ if for every finite
set $\{h_1,\ldots, h_t\}\subseteq \Lambda$
there is a nonzero $a\in A$ such that $a h_i\in A$ for all $i$.
\item[(2)]
We say $A$ is {\it regular pseudo-birational} if, in (1), $a$ can be chosen to be regular in $A$.
\end{enumerate}
\end{definition}

Note that the element $a$ in Definition \ref{xxdef1.3} plays the role of
a common left denominator of
$h_1,\ldots,h_t$. By abuse of language, we 
call $a$ a {\it common left denominator} of
$h_1,\ldots,h_t$.

\begin{example}
\label{xxexa1.4}
\begin{enumerate}
\item[(1)]
Let $I$ be a nonzero right ideal of $\Lambda$
and set $A\colon =\Bbbk+I$. Then $A$ is a
pseudo-birational subalgebra of $\Lambda$.
\item[(2)]
Let $\Lambda=\Bbbk[x]$ and $A$ be the subalgebra generated by $x^2$ and $x^3$. Then $A=\Bbbk+ x^2\Lambda$, and
by part (1), it is regular 
pseudo-birational. However, $A$ is not epimorphic, because
$x\otimes_A 1 -1\otimes_A x$ is a nonzero element in 
the kernel of the map
$\Lambda\otimes_A \Lambda \to \Lambda$.
\end{enumerate}
\end{example}

It is clear that the property of being epimorphic or pseduo-birational is preserved under subalgebra extensions. The next result demonstrates how they behave under taking quotients of algebras.

\begin{lemma}
\label{xxlem1.5}
Let $\Lambda$ be any algebra with some two-sided ideal $N$. Then we have the following. 
\begin{enumerate}
\item[{\rm{(1)}}] 
The image of every {\rm{(}}affine{\rm{)}} epimorphic  subalgebra of $\Lambda$ is again {\rm{(}}affine{\rm{)}} epimorphic in the quotient algebra $\Lambda/N$.
\item[{\rm{(2)}}]
Suppose $A$ is an epimorphic 
subalgebra
of $\Lambda$. Then
$\Lambda\otimes_A(\Lambda/A)=0$.
\item[{\rm{(3)}}] 
Suppose $N$ is finitely generated as a left ideal and it is nilpotent. If $\Lambda/N$ 
admits an affine
regular
pseudo-birational subalgebra, then
$\Lambda$ admits an affine
pseudo-birational subalgebra.
\end{enumerate}
\end{lemma}

\begin{proof}
(1) Let $A$ be an epimorphic subalgebra of $\Lambda$. Denote by $B$ the image of $A$ in the quotient ring $R:=\Lambda/N$. Note that both $A$-actions on
$R$ factor through $B$. Thus we have 
\[
 R\otimes_BR\cong R\otimes_AR\cong R\otimes_\Lambda(\Lambda\otimes_A\Lambda)\otimes_\Lambda R
 \cong R\otimes_\Lambda \Lambda \otimes_\Lambda R  \cong R\otimes_{\Lambda}R\cong R\otimes_{R}R
 \cong R,
\]
which is compatible with the multiplication. Thus $B$ is an epimorphic subalgebra of $R$. 


(2) Since $A$ is epimorphic in $\Lambda$, it follows that $\Lambda\otimes_{A}\Lambda\cong \Lambda$. Moreover, we have that
$$\Lambda \otimes_A(\Lambda/A)
\cong (\Lambda\otimes_A \Lambda)/(\Lambda\otimes_A A)
\cong \Lambda/\Lambda=0.$$

(3) Let $\pi: \Lambda \to \Lambda/N$ be the quotient map. Let $B$ be an affine regular 
pseudo-birational subalgebra of $\Lambda/N$. Choose $b_1,b_2,\dots, b_p
\in \Lambda$ whose images in $\Lambda/N$ generate $B$ as an algebra and $z_1,\dots,z_r\in \Lambda$ generate $N$ as a left ideal. We consider the subalgebra $A$ of $\Lambda$ generated by all $b_i$ and $z_j$. Note that $\pi(A)= B$ and that 
$A$ is affine as well. We establish that $A$ is pseudo-birational in $\Lambda$ by showing that for every finite set $\{h_1,\dots,h_t\}\subseteq\Lambda$, there is a nonzero common left 
denominator $d\in A$ with $\pi(d)\in B$ being regular. Since $N^i\Lambda=N^i$,
we have $$N^{i+1}=\sum_kN^i \Lambda z_k=\sum_k N^i z_k.$$ 
Therefore, the products of length $i$ in the
$z_k$, which also lie in $A$, generate $N^i$ as a left ideal. We construct
$d_i\in A$ such that $\pi(d_i)$ is regular in $B$ and $d_ih_j\in A+N^i$ for all $j$.
For $i=1$, we take a lift $d_1\in A$ of a common left denominator
$b_0\in B$ for the elements $\pi(h_j)\in \Lambda/N$.

Suppose $N^m=0$. For any given $d_i$ with $i<m$, we write $d_ih_j=u_j+n_j$ with $u_j\in A$ and
$n_j=\sum_\ell r_{j\ell}w_\ell\in N^i$, where $r_{j\ell}\in\Lambda$ and the
elements $w_\ell\in A$ are products of length $i$ in the $z_k$. Choose a regular element $b\in B$
with $b\,\pi(r_{j\ell})\in B$ for all $j,\ell$, and a lift $a\in A$ of
$b$. Then $ar_{j\ell}\in A+N$, which implies that $an_j\in A+N^{i+1}$. Hence
$d_{i+1}:=ad_i$ satisfies $d_{i+1}h_j\in A+N^{i+1}$, and
$\pi(d_{i+1})=b\,\pi(d_i)$ is regular in $B$. The claim
follows by taking $d=d_m$. 
\end{proof}

Part (3) of the above lemma is
the third step in \eqref{E0.3.1}, namely,

$$
\text{$\exists$ affine regular pseudo-birational subalgebra in $H/N$}
\Longrightarrow
\text{$\exists$ affine pseudo-birational subalgebra in $H$.}
$$

\section{The main step of the proof}
\label{xxsec2}

In this section, we carry out the main step in \eqref{E0.3.1}, namely,
$$
{\text{$\exists$ affine epimorphic subalgebra in $H$}}
\Rightarrow
{\text{$\exists$ affine regular
pseudo-birational subalgebra  in $H/N$}}.
$$

We start with a baby version of 
this claim. Let $\Ql(A)$ denote 
the classical left quotient
ring of $A$ whenever it exists.

\begin{lemma}
\label{xxlem2.1}
Let $A$ be a subalgebra of $\Lambda$. Suppose that both $A$ and
$\Lambda$ are left Ore domains. Consider the following conditions:
\begin{enumerate}
\item[{\rm{(1)}}]
$A$ is an epimorphic subalgebra.
\item[{\rm{(2)}}]
$\Ql(A)=\Ql(\Lambda)$.
\item[{\rm{(3)}}]
$A$ is a regular pseudo-birational subalgebra.
\end{enumerate}
Then (1) $\Longrightarrow$
(2) $\Longleftrightarrow$ (3).
\end{lemma}

\begin{proof}
It is clear that (2) $\Leftrightarrow$ (3). So we only need to prove that (1) $\Rightarrow$ (2).

Let $E:=\Ql(A)$ and $Q:=\Ql(\Lambda)$.
Then $E$ is a division subalgebra of $Q$, and
$$\begin{aligned}
{\text{Part (1)}}&
\Longrightarrow \Lambda \otimes_A \Lambda \xrightarrow{\cong} \Lambda \qquad\qquad
{\text{by definition}}\\
&\Longrightarrow Q \otimes_A \Lambda \xrightarrow{\cong} Q
\qquad\qquad
{\text{applying $Q\otimes_{\Lambda}-$}}\\
&\Longrightarrow Q \otimes_A Q \xrightarrow{\cong} Q
\qquad\qquad
{\text{applying $-\otimes_{\Lambda} Q$}}\\
&\Longrightarrow m: Q \otimes_E Q \xrightarrow{\cong} Q
\qquad\; 
{\text{because $E=\Ql(A)$}}\\
&
\end{aligned}
$$
where the last map $m$ is the multiplication.
If part (2) failed, there would be an element
$x\in \Lambda\setminus E$. 
Let $\{1, x, y_i\}_i$ be a basis of $Q_E$
and $\{1, x, z_j\}_j$ be a basis of $_EQ$.
Then $x\otimes 1-1\otimes x$ is nonzero in $Q\otimes_E Q$. But 
$x\otimes 1-1\otimes x\in \ker m$,
yielding a contradiction. The assertion
follows.
\end{proof}

However, (2) may not necessarily imply (1) (see Example \ref{xxexa1.4}).
The next subsection deals with the 
non-domain case. 

\subsection{Some classical ring-theoretic results}
\label{xxsec1.2}
Results in the this subsection are 
probably known to experts. Since we cannot 
locate the exact references for some of them, we include their proofs.

We use Goldie's theorem
\cite[Theorem 2.3.6]{MR01}, and the following two classical results.

\begin{lemma}
\label{xxlem2.2}
The following hold.
\begin{enumerate}[label=\textup{(\roman*)}]
\item[{\rm{(1)}}]
\label{FU} 
{\rm{\cite{FU65} and \cite[Corollary
2.3]{Zel73}}}
If $R_0$ is a left order
in a simple artinian ring $Q_0$, there are matrix units $e_{ij}$ of
$Q_0$ and a nonzero non-unital 
subring $T_0$ of their centralizer such that
$T_0e_{ij}\subseteq R_0$ for all $i,j$.
\item[{\rm{(2)}}]
\label{Jat} \cite{Ja69}
A domain over $\kk$
without free pairs is a left and right Ore domain.
\end{enumerate}
\end{lemma}

Let $R$ be a semiprime left Goldie algebra.
Then its classical left quotient ring is 
semisimple artinian, and we write
$\Ql(R)=\prod_{\lambda=1}^vM_{n_\lambda}(D_\lambda)$ where $D_{\lambda}$ are division algebras. The next lemma is a 
slight generalization of Lemma \ref{xxlem2.2}(1).

\begin{lemma}
\label{xxlem2.3}
Retain the notation as above. 
\begin{enumerate}
\item[{\rm{(1)}}]
For each $\lambda$, the matrix units
$e^{(\lambda)}_{ij}$ can be chosen so 
that  there is
$0\neq t_\lambda\in D_\lambda$ with $t_\lambda e^{(\lambda)}_{ij}\in R$
for all $i,j$. 
\item[{\rm{(2)}}]
The division algebra $D_\lambda$ is the centralizer of the matrix units
in the $\lambda$-th factor 
$M_{n_\lambda}(D_\lambda)$.
\end{enumerate}
\end{lemma}

\begin{proof} Part (2) is clear, so we only have to prove part (1).
Let $Q_\lambda$ be the $\lambda$-th factor of $Q:=\Ql(R)$, with identity
$z_\lambda$. Note that $z_{\lambda}$
is central in $Q$. We claim that $R_\lambda=R\cap Q_\lambda$ is a left order
in $Q_\lambda$. Write $z_\lambda=u^{-1}c$ with $u,c\in R$ and $u$ regular; then
$c=uz_\lambda\in R\cap Q_{\lambda}=R_\lambda$. For any given $q\in Q_\lambda$, let us choose $w$
regular in $R$ with $wq\in R$. Then $cw= (uz_{\lambda}) w=uwz_\lambda\in 
R_\lambda$ 
is a
unit of $Q_\lambda$, which implies that
$cw$ is a regular element in $R_{\lambda}$.
Now $(cw)q=c(wq)\in R_\lambda$; 
this proves the claim.
The assertion in part (1) follows by
applying Lemma \ref{xxlem2.2}(1) 
to $(R_{\lambda}, Q_{\lambda})$.
\end{proof}

\begin{lemma}
\label{xxlem2.4}
Let $Q=\prod_\lambda M_{n_\lambda}(D_\lambda)$ with matrix units $e^{(\lambda)}_{ij}$ and
elements $t_\lambda$ as in Lemma~\ref{xxlem2.3}, and let $B$ be a
subalgebra of $Q$ containing $1_Q$ and all $t_\lambda
e^{(\lambda)}_{ij}$. If $B$ does not contain any free pair, then $B$ has a
semisimple artinian classical 
{\rm{(}}left{\rm{)}} 
quotient ring $E\subseteq Q$, and
the regular elements of $B$ are exactly the elements of $B$ that are
units of $Q$.
\end{lemma}

\begin{proof}
Fix a $\lambda$ and, for simplicity, 
write $t=t_{\lambda}$ and 
$e_{ij}=e_{ij}^{(\lambda)}$ for all
$i,j$. Define 
\[
B_{\lambda}
:=\{d\in D_{\lambda}
\mid \ de_{11}\in B\},\qquad T_{\lambda}
:=\kk+B_{\lambda}
\subseteq D_{\lambda}.
\]
Here $\kk$ is identified with $\kk1_{Q_\lambda}$ as well as
$\kk 1_{D_{\lambda}}
\subseteq D_{\lambda}$, since $Q$
is a $\kk$-algebra and scalars commute with the matrix units. Since
$de_{11}\cdot d'e_{11}
=dd'e_{11}$, 
$T_{\lambda}$ is a subalgebra of 
$D_{\lambda}$ and hence a domain. By the definition of
$B$, $t(=t_{\lambda})\in T_{\lambda}$. 
If $x,y\in T_{\lambda}$ 
formed a free pair, then
after subtracting scalars we would be able to find $x,y\in B_{\lambda}$. A relation
$f(xe_{11},ye_{11})=0$ in $B$ gives $f(x,y)e_{11}=0$ after right
multiplication by $e_{11}$; and thus $xe_{11}$ and $ye_{11}$ would form a free pair
in $B$. Hence, by Lemma \ref{xxlem2.2}(2),
$T_{\lambda}$ is an Ore domain; and we denote by $F_{\lambda}
\subseteq D_{\lambda}$ its division ring of fractions.

For $b\in B$ with entries $b_{ij}$ in the $\lambda$-th factor, we have that
$tb_{ij}t
\,e_{11}=(te_{1i})b(te_{j1})\in B$, which implies that $b_{ij}\in
t^{-1}T_{\lambda}
t^{-1}\subseteq F_{\lambda}$. Thus $B\subseteq
E:=\prod_\lambda M_{n_\lambda}(F_\lambda)$.

Next, for $p=\alpha+p_0\in T_{\lambda}$ with $\alpha\in\kk$ and $p_0\in B_{\lambda}$, we have
\begin{equation}
\label{E2.4.1}\tag{E2.4.1}
tpt\,e_{ij}
=\alpha(te_{i1})
(te_{1j})+(te_{i1})
(p_0e_{11})
(te_{1j})
\in B.
\end{equation}

Let $X\in E$ and write $X = \prod_{\lambda}
X_{\lambda}$ 
with $X_{\lambda} = (X^{\lambda}_{ij})
\in M_{n_{\lambda}}(F_{\lambda})$.
Below, we will write them as $X_{ij}$ without the superscript 
$\lambda$.
In each factor, we choose $0\neq r_{\lambda}\in T_{\lambda}$ with
$r_{\lambda}\,tX^{\lambda}_{ij}t^{-1}\in T_{\lambda}$ for all $i,j$. Then $s_\lambda:=(tr_{\lambda}t)1_{Q_\lambda}
=\sum_i tr_{\lambda}t\,e_{ii}$ and $s_\lambda X_\lambda=\sum_{i,j}t(r_{\lambda}tX_{ij}t^{-1})t
\,e_{ij}$ lie in $B$ by \eqref{E2.4.1}. Hence $s:=\sum_\lambda s_\lambda\in B$ is a unit of
$E$ with $sX\in B$.

Note  that $E$ is semisimple artinian. If $b\in B$ is not a unit of $E$, then 
$Xb=0$ for some $0\neq X\in E$ and $(sX)b=0$ with $0\neq sX\in B$. Thus $b$ is not regular 
in $B$. Conversely, units of $E$ in $B$ are regular. A matrix over $F_\lambda$ is invertible over $F_\lambda$ if and only if it is
invertible over $D_\lambda$. Hence $B\cap E^\times=B\cap Q^\times$.
Finally, for $b\in B$ and $c$ regular, writing $bc^{-1}=s^{-1}y$ with
$s,y\in B$ gives $sb=yc$. Therefore, $E=\Ql(B)$.
\end{proof}

Now we can complete the second step in \eqref{E0.3.1}, namely,
$$
{\text{$\exists$ affine epimorphic subalgebra in $H/N$}}\Rightarrow
{\text{$\exists$ affine regular pseudo-birational subalgebra in $H/N$}}
$$

\begin{lemma}
\label{xxlem2.5}
Let $R$ be a semiprime left Goldie algebra that contains no free pairs. Then any epimorphic affine subalgebra of $R$ can be enlarged into an affine regular pseudo-birational subalgebra
$B$ with $\Ql(B)=\Ql(R)$. 
\end{lemma}
\begin{proof}
Let $B$ be an epimorphic affine subalgebra 
of $R$. Without loss of generality, we can 
enlarge $B$ to  contain the finitely many elements $t_\lambda e^{(\lambda)}_{ij}$ 
yielded in Lemma~\ref{xxlem2.3} for $R$. 
Note that $B$ stays affine and 
epimorphic (see the remark before Lemma ~\ref{xxlem1.5}).
Clearly $B$ has no free pairs, so by
Lemma~\ref{xxlem2.4}, $E:=\Ql(B)$ exists, is semisimple and lies in
$Q:=\Ql(R)$. Put $M=R/B$. By Lemma~\ref{xxlem1.5}(2),
\[
0=Q\otimes_R(R\otimes_BM)\cong Q\otimes_E(E\otimes_BM).
\]
Since $E$ is semisimple, $E_E$ is a direct summand of $Q_E$, which implies that 
$E\otimes_BM=0$. Thus every element of $M$ is annihilated by a regular
element of $B$ \cite[\S2.1]{MR01}. As a result, for any $r\in R$, there is a regular
$s\in B$ with $sr\in B$, whence $r=s^{-1}(sr)\in E$. Therefore, $B\subseteq R\subseteq E=\Ql(B)$. Taking the classical left quotient rings, we obtain that $\Ql(B)=\Ql(R)$.

Finally, let $r_1,\dots,r_t\in R$. Since $r_i\in\Ql(B)$, they have a
common left denominator, that is, an element $s\in B$, which is regular in $B$ with $sr_i\in B$ for all $i$. Thus $B$ is regular pseudo-birational in $R$.
\end{proof}

We do not know if Lemma \ref{xxlem2.5} still holds 
when $R$ does contain a free pair. So we leave it as a question 
below.

\begin{question}
\label{xxque2.6}
Suppose $R$ is semiprime left Goldie and contains an affine epimorphic subalgebra $A$. Is there an affine regular pseudo-birational
subalgebra $B\subseteq R$?
\end{question}

In the Hopf algebra setting, we also ask

\begin{question}
\label{xxque2.7}
Let $H$ be a noetherian 
Hopf algebra and $A$ be an affine 
epimorphic subbialgebra. Is $A$ then 
pseudo-birational?
\end{question}

\section{Completion of The proof and Consequences}
\label{xxsec3}

It is well-known that every finite subset of a Hopf algebra $H$ lies in a
finite-dimensional subcoalgebra $C$ of $H$. We denote by $\kk\langle C\rangle$ 
the affine subbialgebra of $H$ generated by $C$. 
For a subalgebra $A\subseteq H$ and a finite-dimensional subcoalgebra $C\subseteq H$, we set
\[
 D_C(A)=\{a\in A\mid \ aC\subseteq A\}.
\]

Recall from Definition \ref{xxdef1.3} that 
a subalgebra $A$ of $H$ is 
{\it pseudo-birational} if 
$D_C(A)\neq 0$ for every finite-dimensional subcoalgebra $C\subseteq H$. 

\subsection{The final step}
\label{xxsec3.1}
The next lemma is the final step in \eqref{E0.3.1}, namely,
$$
{\text{$\exists$ affine pseudo-birational subalgebra in $H$}}\Rightarrow
{\text{$H$ is affine.}}
$$

\begin{lemma}
\label{xxlem3.1}
Let $H$ be a Hopf algebra with bijective antipode $S$.
\begin{enumerate}
\item[(1)]
If $H$ admits an affine pseudo-birational subalgebra, then it admits an affine pseudo-birational subbialgebra.
\item[(2)]
If $A$ is a 
pseudo-birational subbialgebra, then $H=S(A)A$.
\item[(3)]
If $H$ admits an affine pseudo-birational subalgebra,  then $H$ is affine. 
\end{enumerate}
\end{lemma}

\begin{proof}
(1) Let $B$ be an affine pseudo-birational
subalgebra of $H$ generated by a 
finite set $X$. Let $C$ be the finite-dimensional subcoalgebra 
of $H$ generated by $X$. Then 
$A:=\Bbbk\langle C\rangle$ 
is an affine subbialgebra 
of $H$ containing $B$, and 
it is pseudo-birational 
because it contains $B$ 
(a common left denominator 
in $B$ also lies in $A$).

(2)
Fix a finite-dimensional subcoalgebra $C\subseteq H$ and set $D=D_C(A)$. We first show $\Delta(D)\subseteq D\otimes
A$. Let $a\in D$ and $c\in C$. Since $ac_1\in A$ and $A$ is a
subcoalgebra, it follows that
\[
\sum a_1c\otimes a_2
=\sum a_1c_1\otimes a_2c_2S(c_3)
= \sum\Delta(ac_1)\bigl(1\otimes S(c_2)\bigr)
\in (A\otimes A) (1\otimes H) =A\otimes H.
\]
Writing $\Delta(a)=\sum_ip_i\otimes q_i$ with the $q_i\in A$ linearly
independent (as $D\subseteq A$ and $A$ is a coalgebra), we get $p_ic\in A$ for all $i$ and all $c\in C$. Thus $p_i\in D$ by the definition of $D$. 

Secondly we claim that the counit does not vanish on $D$. Since $D\neq 0$ by the definition, we pick any 
$0\neq a\in D$ and write
$\Delta(a)=\sum_ip_i\otimes q_i$ with $p_i\in D$. Then
$a=\sum_i\eps(p_i)q_i$, which implies $\eps(p_i)\neq0$ for some $i$. Now let us choose an element $a\in D$ with $\eps(a)=1$, and write $\Delta(a)=\sum_ip_i\otimes
q_i$ with $p_i\in D$ and $q_i\in A$. Applying the anti-automorphism
$S^{-1}$ to $\sum_i S(p_i)q_i=\epsilon(a)=1$ gives rise to 
$\sum_i S^{-1}(q_i)p_i=1$. Hence
$c=\sum_iS^{-1}(q_i)(p_ic)\in S^{-1}(A)A$ for every $c\in C$. As $C$ is
arbitrary, $H=S^{-1}(A)A$. The assertion follows after applying $S$. 

(3) By part (1), $H$ contains 
an affine pseudo-birational subbialgebra, say $A$, that is generated by a finite set $X$. Now, by part (2), $H$ is generated by $X\cup S(X)$.
\end{proof}

\subsection{Completion of the proof}
\label{xxsec3.2}

\begin{theorem}
\label{xxthm3.2}
Let $H$ be a Hopf algebra with bijective antipode, and let $N$ be its
prime radical. Suppose
\begin{enumerate}
\item
$H^+$ is finitely generated as a left
ideal,
\item
$N$ is nilpotent and finitely generated as a left ideal, and
\item
$H/N$ is left Goldie and does not contain a free pair.
\end{enumerate}
Then $H$ is affine.
\end{theorem}

\begin{proof}
This immediately follows from \eqref{E0.3.1}. 
Note that hypothesis (1) and 
Lemma \ref{xxlem1.2}(2) give an affine
epimorphic subbialgebra of $H$.  It 
passes to the quotient ring $H/N$ by
Lemma \ref{xxlem1.5}(1). Now
hypothesis (3) together with 
Lemma \ref{xxlem2.5} insures 
the existence of an affine regular 
pseudo-birational subalgebra 
$\overline{B}$ of $H/N$.
By hypothesis (2), we can use Lemma
\ref{xxlem1.5}(3) to lift $\overline{B}$
to an affine pseudo-birational 
subalgebra of $H$. Finally, Lemma \ref{xxlem3.1}(3) finishes
the proof.
\end{proof}

\begin{proof}[Proof of Theorem~\ref{xxthm0.2}]
By \cite[Theorem 6.1 and Corollary 6.2]{Sc26}, $S$ is bijective and $H$
is left and right noetherian. Hence both $H^+$ and the prime radical $N$ are
finitely generated left ideals. By Levitzki's theorem \cite[Theorem 2.3.7]{MR01}, $N$ is further nilpotent. Therefore, $H/N$ is semiprime left noetherian and hence left Goldie. A free pair in $H/N$ would lift to a free pair in $H$. As a result, Theorem~\ref{xxthm3.2} applies.
\end{proof}

\subsection{A corollary}
\label{xxsec3.3}

Combining Theorem \ref{xxthm0.2}
with previously known results on
affine noetherian PI Hopf algebras, we obtain the next corollary. Part
(1) answers \cite[Question 5.2]{WZ03}, and part (2) answers a question of
Skryabin \cite{Sk21a} for PI Hopf algebras; see also
\cite[Question 2]{And25}.
We refer to the paper \cite{WZ03}
for all the relevant but unexplained terms.

\begin{corollary}
\label{xxcor3.3}
Suppose $H$ is a noetherian locally PI Hopf algebra. Then $H$ is affine,
and hence PI. Further the following hold.
\begin{enumerate}[label=\textup{(\arabic*)}]
\item $H$ is AS-Gorenstein, Auslander--Gorenstein and
GK--Cohen--Macaulay, and $\GKdim H=\injdim H<\infty$.
\item $H$ is residually finite-dimensional.
\item $H$ is left and right flat over each of its right coideal
subalgebras, and is a projective generator as a left and as a right
module over each of its Hopf subalgebras.
\end{enumerate}
\end{corollary}

\begin{proof}
By Corollary~\ref{xxcor0.3}, $H$ is affine and hence PI. Now, part (1) is
\cite[Theorems 0.1 and 0.2]{WZ03}. Part (2) is then Anan'in's theorem
\cite{Ana92}; see \cite[proof of Corollary 4.6]{Sk21b}. Finally, part (3) follows from \cite[Corollary 0.2]{Sk21b}.
\end{proof}

Brown asked whether every affine PI Hopf algebra is noetherian
\cite[Question C(iii)]{Br07}; and Brown and Gilmartin also asked whether the
two conditions are indeed equivalent for PI Hopf algebras
\cite[Question E]{BGi14}. In a sequel paper \cite{HTWZ}, we will study when
an affine Hopf algebra is noetherian.

\section{Comments, examples and questions}
\label{xxsec4}

In this final section, we provide more examples and give some further comments. We also list a few more questions.

\subsection{Beyond the free-pair hypothesis}
\label{xxsec4.1}


\begin{example}
\label{xxex4.1}
The hypothesis of Theorem~\ref{xxthm0.2} is sufficient but not necessary.
Let $G=\mathbb{Z}^2\rtimes\langle t\rangle$, where $tvt^{-1}=\phi(v)$ for $v\in\mathbb{Z}^2$
and $\phi=\left(\begin{smallmatrix}2&1\\1&1\end{smallmatrix}\right)$.
Then $G$ is polycyclic, and $\Bbbk G$ is thus an affine noetherian Hopf algebra \cite{Ha54}.
Since $\phi$ has the eigenvalue $\lambda=(3+\sqrt5)/2>1$, $G$ is not nilpotent-by-finite. As a result, $G$ has exponential growth 
and contains a free subsemigroup of
rank two.
\end{example}

Thus the no free-pair hypothesis is 
sufficient but not necessary.
Whether it can be dropped is exactly 
Question \ref{xxque0.1}.
The absence of free pair is only used to produce the semisimple
quotient ring in Lemma~\ref{xxlem2.4}. We are wondering if the argument given in the previous section 
can be extended to a more general setting. Below is one piece of positive evidence in this regard.

\begin{proposition}
\label{xxpro4.2}
Let $H$ be a Hopf algebra with bijective antipode, and $N$ a nilpotent
ideal of $H$. Let $A_0$ be an affine subbialgebra with $HA_0^+=H^+$.
Suppose
that an affine subalgebra $T$ with $A_0\subseteq T\subseteq H$ containing a finite generator set of the left ideal $N$ has image
$B$ in $R=H/N$ with a classical left quotient ring $E$, and that $R$
and $E$ embed compatibly in a ring $Q$ such that $E_E$ is a direct
summand of $Q_E$ (for instance, $E$ is semisimple). Then $H$ is affine.
\end{proposition}

\begin{proof}
The surjection $H\otimes_{A_0}H\to H\otimes_TH$ is compatible with
the multiplication, which is bijective on $H\otimes_{A}H$ by
Lemma~\ref{xxlem1.2}(1). Hence the multiplication $H\otimes_TH\to H$ is
bijective and $T$ is an epimorphic
subalgebra. By Lemma \ref{xxlem1.5}(1), $B$ is an epimorphic subalgebra in $R$.
By the proof of Lemma \ref{xxlem2.5}, $B$ is an affine regular pseudo-birational subalgebra of $R$.
Since $T$ is affine with image $B$ and contains a finite generator
set of the left ideal $N$, the proof of Lemma~\ref{xxlem1.5}(3), with $A=T$, shows that
$T$ is pseudo-birational in $H$. The rest of the proof is
\eqref{E0.3.1}.
\end{proof}

\begin{example}
\label{xxex4.3}
Let $H$ be a noetherian Hopf algebra with prime radical $N$ such that
$R:=H/N$ is a domain. Suppose that a subbialgebra $A$ generates $H$ as
a Hopf algebra, i.e., $H$ is generated as an algebra by
$\bigcup_{n\in\mathbb{Z}}S^n(A)$, and that every affine subalgebra of
the image $\overline{A}$ of $A$ in $R$ is contained in an affine left
Ore subalgebra of $R$. This holds, for instance, if $\overline{A}$
contains no free pair by Lemma~\ref{xxlem2.2}(2), or if $A$ is affine
and left noetherian; in both cases $H$ itself may contain free pairs.
Then $H$ is affine.

Indeed, for $a\in A$ we have
$S(a)-\eps(a)=-\sum S(a_1)(a_2-\eps(a_2))\in HA^+$ and
$S^{-1}(a)-\eps(a)=-\sum S^{-1}(a_2)(a_1-\eps(a_1))\in HA^+$. By
induction, $S^n(A)^+\subseteq HA^+$ for all $n\in\mathbb{Z}$, and
hence $H^+=HA^+$. Since $H^+$ is a finitely generated left ideal,
there is an affine subbialgebra $A_0\subseteq A$ with $HA_0^+=H^+$.
The image of $A_0$ in $R$ is contained in an affine left Ore
subalgebra $B$ of $R$. Let $T$ be the subalgebra of $H$ generated by
$A_0$, by lifts of a finite set of generators of $B$, and by a finite
set of generators of the left ideal $N$. Then $T$ is affine with
image $B$ in $R$. Now $E:=\Ql(B)$ is a division subring of the
division ring $Q:=\Ql(R)$, so $E_E$ is a direct summand of $Q_E$, and
Proposition~\ref{xxpro4.2} applies.
\end{example}

\subsection{Related questions}\label{xxsec4.2}

Question \ref{xxque0.1} 
still remains open in the general 
case. Note that affine Hopf algebras 
are not necessarily noetherian. The following is such an example. 

\begin{example}
\label{xxex4.4}
The free algebra $\Bbbk \langle x,y\rangle$ is the universal enveloping algebra of the free Lie
algebra on two generators $x$ and $y$, which is affine but not noetherian.
\end{example}

Therefore, one may ask the following question. 

\begin{question}
\label{xxque4.4}
Let $H$ be an affine 
Hopf algebra with 
$\GKdim H<\infty$. Is 
then $H$ noetherian?
\end{question}

We are also wondering

\begin{question}
\label{xxque4.5}
Can we prove a version of Theorem \ref{xxthm0.2} for other classes of Hopf-like algebras (e.g.
weak Hopf algebras, quasi-Hopf algebras, 
Hopf algebras in braided categories, 
Hopf superalgebras and Hopf algebroids)?
\end{question}

We refer the reader to a most recent preprint \cite{JZ26c} for some progress on the above question concerning connected braided Hopf algebras. Other related questions can be found 
in the survey papers mentioned in the 
introduction.

\subsection*{Acknowledgments}

This project was motivated by a talk given 
by Ken Brown at the Seattle Noncommutative Algebra Conference
during December 15--19, 2025. We thank him for the inspiration, and the organizers of the conference for
the opportunity to  
discuss many open problems.
The authors acknowledge the 
use of AI tools in gathering initial 
ideas. The final proofs were 
independently developed and 
checked by the authors, 
who take full responsibility 
for the integrity of the 
mathematical content. X.W.\ was partly supported by Simons Travel Support for Mathematicians SFI-MPS-TSM-00025285.
James J. Zhang was partially supported by the US National Science Foundation (No. DMS-2302087).

\end{document}